\documentclass[12pt]{article}
\usepackage{amssymb,amsfonts,amsmath, psfrag}
\newtheorem{theo}{Theorem}

\newtheorem{lem}[theo]{Lemma}

\newtheorem{pro}[theo]{Proposition}

\def\qed{\hfill \rule{4pt}{7pt}}
\def\MX{{\mathcal{X}}}
\def\MC{{\mathcal{C}}}
\def\MB{{\mathcal{B}}}
\begin{document}
\begin{center}
{\bf \Large Rational Interpolation and Basic Hypergeometric
Series}
\end{center}
\vspace{0.3cm}
\begin{center}
  {\bf Amy M. Fu}\\
  Center for Combinatorics\\
  Nankai University, Tianjin 300071, P.R. China\\
  Email: fmu@eyou.com
 \end{center}
\begin{center}
 {\bf Alain Lascoux}\\
 Nankai University, Tianjin 300071, P.R. China\\
 Email: Alain.Lascoux@univ-mlv.fr\\
 CNRS, IGM Universit\'e  de
 Marne-la-Vall\'ee\\
 77454 Marne-la-Vall\'ee Cedex, France\\
 \end{center}

 \vspace{1cm}
 \noindent{\bf Abstract.} We give a Newton type rational interpolation formula (Theorem \ref{theo}).
 It contains as a special case the original Newton interpolation, as
 well as the recent interpolation formula of Zhi-Guo Liu, which
 allows to recover many important classical $q$-series identities.
 We show in particular that some bibasic identities are a
 consequence of our formula.

\section{Introduction and Notation}

As usual, $(a;q)_n$ denotes:
 $$\prod_{j=0}^{n-1}(1-aq^{j}),
n=0,1,2,\cdots,\infty,$$ and
$$
(a_1,a_2,\cdots,a_m;q)_n=(a_1;q)_n(a_2;q)_n\cdots(a_m;q)_n.
$$
The basic hypergeometric series $_{r}\phi_s$ are defined in
\cite{Gas} by:
\begin{eqnarray*}
_r\phi_s\left[\begin{array}{c} a_1,\ a_2, \ \cdots, \ a_r\\
b_1,\ b_2, \ \cdots, \ b_s
\end{array};q,
z\right]=\sum_{n=0}^{\infty}\frac{(a_1,a_2,\cdots,a_r;q)_n}{(b_1,b_2,\cdots,b_s;q)_n}
\left[(-1)^nq^{\binom{n}{2}}\right]^{1+s-r}z^n,
\end{eqnarray*}
where $q\neq 0$ when $r>s+1$.

Newton obtained the following interpolation formula:
\begin{eqnarray*}
f(x)=f(x_1)+f^{\partial}(x-x_1)+f^{\partial\partial}(x-x_1)(x-x_2)+\cdots,
\end{eqnarray*}
where $f^{\partial\cdots \partial}$ is the image of $f(x)$ under a
product of divided differences which will be defined below.

Special cases of Newton's interpolation are the Taylor formula and
the $q$-Taylor formula (c.f. \cite{Kac}), with derivatives or
$q$-derivatives instead of divided differences.

Using $q$-derivatives, Zhi-Guo Liu \cite{Liu} gave  an
interpolation formula involving rational functions in $x$ as
coefficients, instead of only polynomials in $x$ as in the
$q$-Taylor formula:
\begin{equation}\label{liu}
f(x)=\sum_{n=0}^{\infty}\frac{(1-aq^{2n})(aq/x;q)_{n}x^n}{(q,x;q)_{n}}[D_{q}f(x)(x;q)_{n-1}]{\big|}_{x=aq},
\end{equation}
$D_q$ being defined by
$$
D_qf(x)=\frac{f(x)-f(xq)}{x}.
$$

Let us remark that Carlitz's $q$-analog of a special case of the
Lagrange inversion formula is the limit for $a \rightarrow 0$ of
(\ref{liu}):
\begin{eqnarray*}
f(x)=\sum_{n=0}^{\infty}\frac{x^n}{(q,x;q)_{n}}[D_{q}f(x)(x;q)_{n-1}]|_{x=0}.
\end{eqnarray*}

Our formula involves two sets of indeterminate ${\mathcal X}$ and
${\mathcal C}$. Newton interpolation is the case when
$${\mathcal C}=\{0,0,\cdots\},$$
 and Zhi-Guo Liu's expansion is the
case when
$$
{\mathcal X}=\{aq^1, aq^2,\cdots \}, \ {\mathcal
C}=\{q^0,q^1,q^2,\cdots\}.
$$

\section{Rational Interpolation}

By convenience, we denote
$$
Y_n(x,\MX)=(x-x_1)(x-x_2)\cdots(x-x_n)$$
and
$$
(x,\MC)_n=(1-xc_1)(1-xc_2)\cdots(1-xc_n).
$$

The divided difference $\partial_i$ (acting on its left),
$i=1,2,3,\cdots$ is defined by
\begin{multline*}
f(x_{1},x_{2},\cdots, x_{i},x_{i+1},\cdots)\,\partial_{i}
\\
=\frac{f(\cdots, x_{i},x_{i+1},\cdots)-f(\cdots,
x_{i+1},x_{i},\cdots)}{x_{i}-x_{i+1}}=\frac{f-f^{s_i}}{x_i-x_{i+1}},
\end{multline*}
where $s_i$ denotes the exchange of $x_i$ and $x_{i+1}$.

 Divided differences satisfy  a Leibnitz type formula:
 \begin{eqnarray*}
(f(x_1)g(x_1))\partial_1= f(x_1)\, \bigl( g(x_1)\partial_1 \bigr)+
\bigl(f(x_1)\partial_1 \bigr)\, g(x_2).
\end{eqnarray*}
By induction, one obtains:
\begin{eqnarray*}
f(x_1)g(x_1)\,
\partial_1\partial_2\cdots\partial_n=\sum_{k=0}^n\bigl(
f(x_{1})\partial_{1}\cdots\partial_k\bigr)\,\bigl(g(x_{k+1})\partial_{k+1}\cdots\partial_n\bigr).
\end{eqnarray*}

\begin{lem}\label{lem1}Letting $i$,$n$ be two nonnegative integers, then one has:
$$
Y_n(b_1,\MX) \,
\partial_1\partial_2\cdots\partial_i{\big|}_{\MB=\MX}=\left\{\begin{array}{l}
0, \ i \neq n;\\
1, \  i=n,
\end{array}\right.
$$
where ${\big |}_{\MB=\MX}$ denotes the specialization $b_1=x_1,\
b_2=x_2, \cdots,$ and the divided differences are relative to
$b_1, b_2, \cdots$.
\end{lem}

\noindent{\bf Proof.} If $i\leq n$, using Leibnitz formula, we
have:
\begin{eqnarray*}
&&Y_n(b_1,\MX) \,\partial_1\partial_2\cdots\partial_i{\big|}_{\MB=\MX}\\
&&=\prod_{k=2}^n(b_1-x_k) \,
\partial_1\partial_2\cdots\partial_i
(b_1-x_1){\big|}_{\MB=\MX}+\prod_{k=2}^n(b_2-x_k) \,
\partial_2\cdots\partial_i
(b_1-x_1) \,\partial_1{\big|}_{\MB=\MX}\\
&&=\prod_{k=2}^n(b_2-x_k) \, \partial_2\cdots\partial_i{\big|}_{\MB=\MX}
= \cdots =\prod_{k=i}^n(b_{i}-x_k) \, \partial_i{\big|}_{\MB=\MX}\\
&&=\left\{\begin{array}{l}
\prod_{k={i+1}}^n(b_{i+1}-x_k){\big|}_{\MB=\MX}=0, \ i<n;\\
(b_n-x_n) \, \partial_n{\big|}_{\MB=\MX}=1, \  i=n.
\end{array}\right.
\end{eqnarray*}

In the case $i>n$, nullity comes from the fact that each
$\partial_i$ decreases degree by $1$. \qed

\begin{theo}\label{theo}
For any formal series $f(x)$ in $x$, we have the following
identity in the ring of formal series in $x, x_1, x_2, \cdots$:
\begin{eqnarray}\label{las1}
f(x)&=&f(x_{1})+f(x_1) \, \partial_{1}(1-x_{2}c_{1})\frac{Y_{1}(x,\MX)}{(x,\MC)_{1}}\nonumber\\
&+&f(x_1)(1-x_{1}c_{1})\,
\partial_{1}\partial_{2}(1-x_{3}c_{2})\frac{Y_{2}(x,\MX)}{(x,\MC)_2}\nonumber\\
&+&\cdots + f(x_1)(x_{1},\MC)_{n-1} \,
\partial_{1}\cdots\partial_{n}(1-x_{n+1}c_n)\frac{Y_{n}(x,\MX)}{(x,\MC)_{n}}+\cdots.
\end{eqnarray}
\end{theo}

\noindent{\bf Proof.} Let
$$
f(b)=\sum_{n=0}^{\infty}A_n\frac{Y_{n}(b,\MX)}{(b,\MC)_{n}}.
$$
Specializing $b$ to $x_1$ or $x_2$, one gets the following
coefficients:
$$
A_0=f(x_1), A_1=f(x_1)\, \partial_1 (1-x_2c_1).
$$

 Now we have to check,
\begin{eqnarray}\label{Mid}
\frac{Y_n(b_1,\MX)}{(b_1,\MC)_n}(b_1,\MC)_{k-1} \,
\partial_1\partial_2\cdots\partial_k
{\big|}_{\MB=\MX}=\left \{\begin{array}{cc}0, & k\neq n;\\
\frac{1}{1-x_{n+1}c_n}, & k=n.\end{array} \right.
\end{eqnarray}

If $k>n$,
$\displaystyle{\frac{Y_n(b,\MX)}{(b,\MC)_n}(b,\MC)_{k-1}}$ is a
polynomial of degree $k-1$, and therefore annihilated by a product
of $k$ divided differences.

If $k<n$, from Leibnitz formula, we get:
\begin{eqnarray*}
&&\frac{Y_n(b_1,\MX)}{(b_1,\MC)_n}(b_1,\MC)_{k-1} \,
\partial_1\partial_2\cdots\partial_k
{\big|}_{\MB=\MX}\\
&&=\frac{Y_n(b_1,\MX)}{\prod_{p=k}^n(1-b_1c_p)} \,
\partial_1\partial_2\cdots\partial_k
{\big|}_{\MB=\MX}\\
&&=\sum_{i=0}^k\frac{1}{\prod_{p=k}^n(1-b_{i+1}c_p)} \,
\partial_{i+1}\cdots
\partial_k Y_n(b_1,\MX) \,
\partial_1\cdots\partial_i{\big|}_{\MB=\MX},
\end{eqnarray*}
and Lemma \ref{lem1} shows that this function is equal to $0$.

If $k=n$, we have:
\begin{eqnarray*}
&&\frac{Y_n(b_1,\MX)}{1-b_1c_n} \,
\partial_1\partial_2\cdots\partial_n
{\big|}_{\MB=\MX}\\
&&=\sum_{i=0}^n\frac{1}{1-b_{i+1}c_n} \, \partial_{i+1}\cdots
\partial_n Y_n(b_1,\MX) \, \partial_1\cdots\partial_i{\big|}_{\MB=\MX}\\
&&=\frac{1}{1-b_{n+1}c_n}Y_n(b_1,\MX)\, \partial_1\cdots\partial_n{\big|}_{\MB=\MX}\\
&&=\frac{1}{1-x_{n+1}c_n}.
\end{eqnarray*}
Formula (\ref{Mid}) thus implies
$$
A_n=f(x_1)(x_1,\MC)_{n-1}\partial_1\cdots\partial_n(1-x_{n+1}c_n),
$$
and the theorem. \qed

\section{Application to Basic summation formulas}

Taking
 $$f(x)=(x;q)_n,$$
 and
 $$\MX=\{aq^1,aq^2,\cdots\},\, \MC=\{0,0,\cdots\},$$
we have:
\begin{eqnarray*}
(x;q)_n&=&\sum_{k=0}^n(x-aq)\cdots(x-aq^k)(x_1;q)_n\partial_1
\cdots \partial_k \\
&=&\sum_{k=0}^n(x-aq)\cdots(x-aq^{k})(-1)^kq^{\binom{k}{2}}{n
\brack k}
(aq^{k+1};q)_{n-k}\\
&=&(aq;q)_n\sum_{k=0}^n\frac{(q^{-n},aq/x;q)_k}{(q,aq;q)_k}(xq^n)^k.
\end{eqnarray*}
Replacing $aq$ by $a$, $x$ by $a/c$,  we get a $q$-analogue of
Vandermonde's formula:
\begin{eqnarray*}
\frac{(a/c;q)_n}{(a;q)_n}=\sum_{k=0}^n\frac{(q^{-n},c;q)_k}{(q,a;q)_k}(aq^n/c)^k=
{_2}\phi_1\left[\begin{array}{cc} q^{-n}, &c\\
a
\end{array};q,aq^n/c
\right].
\end{eqnarray*}
Recall another $q$-analogue of Vandermonde's formula:
\begin{equation}\label{Van}
\frac{(a/c;q)_n}{(a;q)_n}c^n=\sum_{k=0}^n\frac{(q^{-n},c;q)_k}{(q,a;q)_k}q^k.
\end{equation}

 Taking
$$f(x)=\frac{1}{(a\beta x;q)_{\infty}}
$$
and
$$\MX=\{q/a, q^2/a, \cdots \}, \MC=\{a\beta, a\beta q,
\cdots\},$$ gives:
$$
f(x_1)(x_1;q)_{n-1}\partial_{1}\cdots\partial_{n}=\frac{1}{(a\beta
x_1q^{n-1};q)_{\infty}}\partial_{1}\cdots\partial_{n}
=\frac{(a\beta)^nq^{n(n-1)}}{(q;q)_{n}(\beta q^{n};q)_{\infty}}.
$$

From Theorem \ref{theo}, we have:
\begin{eqnarray*}
\frac{1}{(a\beta
x;q)_{\infty}}=\frac{1}{(\beta;q)_{\infty}}\sum_{n=0}^{\infty}\frac{(1-\beta
q^{2n})(q/ax,\beta;q)_{n}(a\beta x)^nq^{n(n-1)}}{(q;q)_{n}(a\beta
x;q)_{n}},
\end{eqnarray*}
and we get Jackson's formula \cite{Jack}:
\begin{equation}\label{la6}
\frac{(\beta;q)_{\infty}}{(a\beta
x;q)_{\infty}}=\sum_{n=0}^{\infty}\frac{(1-\beta q^{2n})(q/ax,
\beta; q)_{n} (a\beta x q^{n-1})^n}{(q;q)_{n}(a\beta x;q)_{n}}.
\end{equation}
Setting $x \rightarrow 0$, replacing $\beta$ by $\beta q$, we
derive Sylvester's formula:
\begin{eqnarray*}
\sum_{n=0}^{\infty}\frac{(-1)^n\beta^n q^{n(3n+1)/2}(1-\beta
q^{2n+1} )}{(q;q)_{n}(\beta q^{n+1};q)_{\infty}}=1.
\end{eqnarray*}

In (\ref{la6}), letting $a=1$, $x=q^{N}$, replacing $\beta$ by
$-\beta q$, we get Andrews' formula \cite{Andr2}:
\begin{eqnarray*}
(-\beta q;q)_{2N}=\sum_{n=0}^{N}(-\beta q;q)_{n-1}(1+\beta
q^{2n})q^{n(3n-1)/2}{N \brack n}(-\beta q^{n+N+1};q)_{N-n}.
\end{eqnarray*}

 Taking
$$
f(x)={_3}\phi_2\left[\begin{array}{ccc} q^{-n}, &a, &x/b\\
x, &aq^{1-n}/e &
\end{array};q,q
\right],
$$
and
 $$\MX=\{0,0,\cdots\},\  \MC=\{q^0,q^1,q^2,\cdots\},$$

we have:
\begin{eqnarray*}
f(x)&=&\sum_{i=0}^{\infty}\frac{x^i}{(x;q)_i}f(x_1)(x_1;q)_{i-1} \,\partial_1\cdots\partial_i\\
&=&\sum_{i=0}^{\infty}\frac{x^i}{(x;q)_i}
\sum_{k=i}^{n}\frac{(q^{-n},a;q)_kq^k}{(q,aq^{1-n}/e;q)_k}
 \frac{(x_1/b;q)_k}{(x_1q^{i-1};q)_{k-i+1}} \,\partial_1\cdots\partial_i\\
&=&\sum_{i=0}^{\infty}\frac{(q^{-n},a,q^{1-i}/b;q)_i(xq^{i})^i}{(q,x,aq^{1-n}/e;q)_i}
\sum_{k=i}^{n}\frac{(q^{-n+i},aq^i;q)_{k-i}q^{k-i}}{(q,aq^{1-n+i}/e;q)_{k-i}}\\
&=&\sum_{i=0}^{\infty}\frac{(q^{-n},a,q^{1-i}/b;q)_i(xq^{i})^i}{(q,x,aq^{1-n}/e;q)_i}
\frac{(q^{1-n}/e;q)_{n-i}}{(a/eq^{1-n+i};q)_{n-i}}(aq^i)^{n-i} \ \ ({\rm from} \ \ref{Van})\\
&=&\frac{(q^{1-n}/e;q)_{n}}{(a/eq^{1-n};q)_{n}}a^n\sum_{i=0}^{n}\frac{(q^{-n},
a,b;q)_i}{(q,x,e;q)_i}\left(xeq^n\over ab\right)^i.
\end{eqnarray*}
Raplacing $x$ by $d$, we obtain Sears' formula (\cite{Gas}
(3.2.5)):
\begin{eqnarray*}
{_3}\phi_2\left[\begin{array}{ccc} q^{-n}, &a, &b\\
d, &e &
\end{array};q,deq^n/ab
\right]=\frac{(e/a;q)_n}{(e;q)_n}{_3}\phi_2\left[\begin{array}{ccc} q^{-n}, &a, &d/b\\
d, &aq^{1-n}/e &
\end{array};q,q
\right].
\end{eqnarray*}

\section{Bibasic summation formula}

\begin{pro}\label{pro}
Taking
$$
f(x)=\frac{1-ux}{1-vx}
$$
and
$$
\MX=\{x_1, x_2, \cdots\}, \ \MC=\{c_1, c_2, \cdots\},
$$
 we have:
\begin{equation}\label{Pro}
f(x)=\frac{1-ux_1}{1-vx_1}+\sum_{k=1}^{\infty}\frac{(v-u)Y_{k-1}(v,\MC)(1-x_{k+1}c_k)}{(v,\MX)_{k+1}}
\frac{Y_{k}(x,\MX)}{(x,\MC)_k}
\end{equation}
\end{pro}

The proposition is a direct application of Theorem \ref{theo} and
of the following lemma.

\begin{lem}\label{Lem}
\begin{equation}\label{Main}
\frac{1-ux_1}{1-vx_1}(x_1,\MC)_{k-1}\partial_1\partial_2\cdots\partial_k=(v-u)Y_{k-1}(v,\MC)_k/(v,\MX)_{k+1}.
\end{equation}
\end{lem}

We first need to recall some facts about symmetric functions
\cite{M}. The generating function for the elementary symmetric
functions $e_i(x_1, x_2, \cdots)$, and the complete functions
$S_i(x_1, x_2, \cdots)$ are
\begin{eqnarray*}
\sum_{i\geq 0} e_i(x_1, x_2, \cdots)t^i=\prod_{i \geq 0}(1+x_it),
\end{eqnarray*}
and
\begin{eqnarray*}
\sum_{i\geq 0} S_i(x_1, x_2, \cdots)t^i=\prod_{i \geq
0}1/(1-x_it).
\end{eqnarray*}

We shall need a slightly more general notion than usual, for a
Schur function. Given $\lambda \in \mathbb{N}^n$, given $n$ sets
of variable $A_1, \ldots, A_n$, then the Schur function
$S_{\lambda}(A_1, \ldots, A_n)$ is equal to $
\bigl|S_{\lambda_j+j-i}(A_j)\bigr|_{1\leq i,j\leq n}. $

One has the following identity \cite{Lascoux}:
\begin{equation}\label{Las}
S_{\lambda}(x_2, x_3, \cdots)x_1^r=S_{\lambda, r}(\MX,x_1).
\end{equation}

 \noindent{\bf Proof of Lemma \ref{Lem}} Multiply the
denominator of $\displaystyle{\frac{1-ux_1}{1-vx_1}(x_1,
\MC)_{k-1}}$ by the symmetrical factor $(v,\MX)_{k+1}$, which
commutes with $\partial_1\cdots
\partial_k$. Let ${\MX}_k=\{x_1, x_2,\ldots x_{k+1}\}$. One has:
\begin{eqnarray*}
&&(1-ux_1)\prod_{i=1}^{k-1}(1-x_1c_i)
\prod_{j=2}^{k+1}(1-vx_j)\\
&&=\sum_{i=0}^{k}\sum_{j=0}^k(-1)^i(-v)^je_i(u,c_1,\cdots,c_{k-1})e_j(x_2,x_3,\cdots,
x_{k+1})x_1^i\\
&&=\sum_{i=0}^{k}\sum_{j=0}^k(-1)^i(-v)^je_i(u,c_1,\cdots,c_{k-1})S_{1^j,i}(\underbrace{{\MX}_k,\ldots,{\MX}_k}_j,
x_1),
\end{eqnarray*}
thanks to (\ref{Las}), and to the fact that for every $j$,
$e_j(\MX)=S_{1^j}(\underbrace{\MX,\ldots,\MX}_j)$.

The image of a power of $x_1$ under $\partial_1\cdots\partial_k$
is a complete function in $\MX$ \cite{Lascoux}. Therefore,
$$
S_{1^j,i}(\underbrace{{\MX}_k,\ldots,{\MX}_k}_j,x_1) \,\partial_1
\cdots
\partial_k=S_{1^j,i-k}(\underbrace{{\MX}_k,\ldots,{\MX}_k}_{j+1}).
$$
This determinant is equal to $0$ (because it has two identical
columns), except for $i+j=k$, in which case it is equal to
$S_{0^{j+1}}(\MX)=(-1)^j$.

Now:
\begin{multline*}
\frac{1-ux_1}{1-vx_1}(x_1,\MC)_{k-1}(v,\MX)_{k+1} \,\partial_1\partial_2\cdots\partial_k\\
=\sum_{i+j=n}(-1)^iv^je_i(u,c_1,\cdots,c_{k-1})=(v-u)Y_{k-1}(v,\MC)_k,
\end{multline*}
thus (\ref{Main})  is true. \qed

In \cite{GasF}, Gasper obtained the following identity:
\begin{equation}\label{Gas}
\sum_{k=0}^{\infty}\frac{1-ap^kq^k}{1-a}\frac{(a;p)_k(b^{-1};q)_kb^k}{(q;q)_k(abp;p)_k}=0.
\end{equation}

We also prove an identity due to Gosper (c.f. \cite{Gas}):
\begin{eqnarray*}
\sum_{k=0}^n\frac{1-ap^kq^k}{1-a}\frac{(a;p)_k(c;q)_kc^{-k}}{(q;q)_k,(ap/c;p)_k}=
\frac{(ap;p)_n(cq;q)_nc^{-n}}{(q;q)_n(ap/c;p)_n},
\end{eqnarray*}
or equivalently,
\begin{equation}\label{Gos}
\sum_{k=0}^n\frac{(1-ap^{n-k}q^{n-k})(q^{n-k+1};q)_{k}(ap^{n-k+1}/c;p)_{k}}
{(cq^{n-k};q)_{k+1}(ap^{n-k};p)_{k+1}}c^{k}=\frac{1}{1-c}.
\end{equation}

In fact, (\ref{Gas}) and (\ref{Gos}) are  special cases of
Proposition \ref{pro}.

Taking $u=0$ in (\ref{Pro}), we get:
\begin{equation}\label{Gas1}
\frac{1}{1-vx}=\frac{1}{1-vx_1}+\sum_{k=1}^{\infty}\frac{vY_{k-1}(v,\MC)(1-x_{k+1}c_k)}{(v,\MX)_{k+1}}
\frac{Y_{k}(x,\MX)}{(x,\MC)_k}.
\end{equation}
 Multiplying both sides of (\ref{Gas1}) by $(1-vx_1)$, one has:
\begin{eqnarray*}
\frac{1-vx_1}{1-vx}=1+\sum_{k=1}^{\infty}\frac{vY_{k-1}(v,\MC)(1-x_{k+1}c_k)}{(v,\MX/x_1)_{k}}
\frac{Y_{k}(x,\MX)}{(x,\MC)_k},
\end{eqnarray*}
where $\MX/x_1=\{x_2,x_3,\cdots\}$.

Taking
$$\MX=\{q^0,q^1,q^2,\cdots\},\  \MC=\{ap^1, ap^2,\cdots\}, \ v=1, \ x=b,$$ we get (\ref{Gas}).

In (\ref{Pro}), taking
$$
\MX=\{p^{-n}/a, p^{-n+1}/a,\cdots\}, \
\MC=\{q^{-n+1},q^{-n+2},\cdots\},
$$
$u=q^{-n}$, $v=1$, $x=c^{-1}$, we get (\ref{Gos}).

 \vskip 2mm \noindent{\bf Acknowledgments.}

This work was done under the auspices of the National Science
Foundation of China.

\end{document}